# On the resolution and Linear programming problems subjected by Aczel-Alsina Fuzzy relational equations

**Amin Ghodousian*, Hadi Amiri, Alireza Norouzi Azad**

Faculty of Engineering Science, College of Engineering, University of Tehran, P.O.Box 11365-4563, Tehran, Iran ;
*a.ghodousian@ut.ac.ir*

Faculty of Engineering Science, College of Engineering, University of Tehran, P.O.Box 11365-4563, Tehran, Iran ;
*hadi.amiri@ut.ac.ir*

Department of Engineering Science, College of Engineering, University of Tehran, Tehran, Iran;
*alireza.norouzi@ut.ac.ir*

**Abstract**

Aczel-Alsina t-norm belongs to the family of strict t-norms that are the most applied fuzzy operators in various fuzzy modelling problems. In this paper, we study a linear optimization problem where the feasible region is formed as a system of fuzzy relational equations (FRE) defined by the Aczel-Alsina t-norm. Since the feasible solutions set of FREs is non-convex and the finding of all minimal solutions is an NP-hard problem, conventional methods may not be directly employed. The resolution of the feasible region is completely investigated. Based on some theoretical properties of the problem, an algorithm is presented to find all the optimal solutions, and finally an example is described to illustrate this algorithm.

## 1. Introduction

In this paper, we study the following linear problem in which the constraints are formed as fuzzy relational equations defined by Aczel-Alsina t-norm:

$$\min \ cx$$
$$\varphi(a_i, x) = b_i \ , \ i \in I \tag{1}$$
$$x \in [0,1]^n$$

where $I = \{1,2,...,m\}$ , $J = \{1,2,...,n\}$, $0 \leq a_{ij} \leq 1$ ($\forall i \in I$ and $\forall j \in J$), $b = (b_i)_{m \times 1}$, $0 \leq b_i \leq 1$ ($\forall i \in I$), "$\varphi$" is the max-Aczel-Alsina composition and the constraints mean:

$$\varphi(a_i, x) = \max_{j \in J}\{\varphi(a_{ij}, x_j)\} = \max_{j \in J}\{T_{AA}^{\lambda}(a_{ij}, x_j)\} = b_i \quad , \forall i \in I$$

and

---

* Corresponding author
Email addresses: a.ghodousian@ut.ac.ir (Amin Ghodousian), Hadi.Amiri@ut.ac.ir (Hadi Amiri), alireza.Norouzi@ut.ac.ir (Alireza Norouzi Azad).

$$T_{AA}^{\lambda}(a_{ij},x_j) = \begin{cases} 0 & , a_{ij}=0 \text{ or } x_j=0 \\ \exp\left(-\left[\left(-Ln\,(a_{ij})\right)^{\lambda}+\left(-Ln\,(x_j)\right)^{\lambda}\right]^{\frac{1}{\lambda}}\right) & \text{otherwise} \end{cases}$$

It can be easily shown that Aczel-Alsina t-norm $T_{AA}^{\lambda}(x,y)$ converges to the minimum fuzzy intersection $\min\{x,y\}$ as $\lambda$ goes to infinity and converges to Drastic product t-norm as $\lambda$ approaches $0$ [2]. Also, it is interesting to note that $T_{AA}^{1}(x,y)=xy$, that is, the Aczel-Alsina t-norm is converted to the product t-norm if $\lambda=1$.

The theory of fuzzy relational equations (FRE) was firstly proposed by Sanchez and applied in problems of the medical diagnosis [13]. Nowadays, it is well known that many issues associated with a body knowledge can be treated as FRE problems [12]. The solution set of FRE is often a non-convex set that is completely determined by one maximum solution and a finite number of minimal solutions [4]. The other bottleneck is concerned with detecting the minimal solutions that is an NP-hard problem [5,7,8,10]. The problem of optimization subject to FRE and FRI is one of the most interesting and on-going research topic among the problems related to FRE and FRI theory [1,2,4–8,11,14]. Recently, many interesting generalizations of the linear programming subject to a system of fuzzy relations have been introduced and developed based on composite operations used in FRE, fuzzy relations used in the definition of the constraints, some developments on the objective function of the problems and other ideas [3,7-9,11].

In this paper, an algorithm is proposed to find all the optimal solutions of problem (1). For this purpose, we describe some structural details of FREs defined by the Aczel-Alsina t-norm such as the theoretical properties of the solutions set, necessary and sufficient conditions for the feasibility of the problem and some simplification processes to reduce the problem.

The remainder of the paper is organized as follows. Section 2 gives some basic results on the feasible solutions set of problem (1). The optimal solution of the problem is characterized in Section 3 and, finally in section 4 the experimental results are demonstrated.

## 2. Feasible region of max-Aczel-Alsina FRE

Let $S(a_{ij},b_i)=\{x_j\in[0,1]:T_{AA}^{\lambda}(a_{ij},x_j)=b_i\}$ and $S(a_i,b_i)$ denote the feasible solutions set of $i$'th equation, that is, $S(a_i,b_i)=\{x\in[0,1]^n:\max_{j=1}^{n}\{T_{AA}^{\lambda}(a_{ij},x_j)\}=b_i\}$. Also, let $S(A,b)$ denote the feasible solutions set of problem (1). Based on the foregoing notations, it is clear that $S(A,b)=\bigcap_{i\in I}S(a_i,b_i)$. Also, for each $i\in I$, we define $J_i=\{j\in J:a_{ij}\geq b_i\}$.

Obviously, if $a_{ij}=b_i=0$, then $T_{AA}^{\lambda}(a_{ij},x_j)=b_i$ for each $x_j\in[0,1]$. Also, if $a_{ij}>b_i=0$, then $S(a_{ij},b_i)=\{0\}$; because $0$ satisfies the equality $T_{AA}^{\lambda}(a_{ij},x_j)=b_i$ and $T_{AA}^{\lambda}$ is a strict t-norm. Moreover, from the definition of Aczel-Alsina t-norm, if $a_{ij}>b_i>0$, then $\exp(-[(-Ln\,(b_i))^{\lambda}-(-Ln\,(a_{ij}))^{\lambda}]^{1/\lambda})$ is the unique solution to the equality $T_{AA}^{\lambda}(a_{ij},x_j)=b_i$. Particularly, if $a_{ij}=b_i>0$, we have $S(a_{ij},b_i)=\{1\}$. The foregoing results have been summarized in the following lemma.

**Lemma 1.** Let $i \in I$. Then, $S(a_{ij}, b_i) = \emptyset$ if $j \notin J_i$; $S(a_{ij}, b_i) = [0,1]$ if $a_{ij} = b_i = 0$; $S(a_{ij}, b_i) = \{1\}$ if $a_{ij} = b_i > 0$; $S(a_{ij}, b_i) = \{0\}$ if $a_{ij} > b_i = 0$; $S(a_{ij}, b_i) = \{\exp(-[(-Ln(b_i))^\lambda - (-Ln(a_{ij}))^\lambda]^{1/\lambda})\}$ if $a_{ij} > b_i > 0$.

**Definition 1.** Suppose that $i \in I$ and $S(a_i, b_i) \neq \emptyset$. Let $\hat{x}_i = [(\hat{x}_i)_1, (\hat{x}_i)_2, ..., (\hat{x}_i)_n] \in [0,1]^n$ where the components ($\forall k \in J$) are defined as follows:

$$(\hat{x}_i)_k = \begin{cases} \exp(-[(-Ln(b_i))^\lambda - (-Ln(a_{ik}))^\lambda]^{1/\lambda}) & , k \in J_i, b_i \neq 0 \\ 0 & , k \in J_i, a_{ik} > b_i = 0 \\ 1 & , k \in J_i, a_{ik} = b_i = 0 \end{cases}$$

Also, for each $j \in J_i$, we define $\breve{x}_i(j) = [\breve{x}_i(j)_1, \breve{x}_i(j)_2, ..., \breve{x}_i(j)_n] \in [0,1]^n$ such that for each $k \in J$:

$$\breve{x}_i(j)_k = \begin{cases} \exp(-[(-Ln(b_i))^\lambda - (-Ln(a_{ik}))^\lambda]^{1/\lambda}) & , b_i \neq 0 \text{ and } k = j \\ 0 & \text{otherwise} \end{cases}$$

**Definition 2.** Let $\hat{x}_i$ ($i \in I$) be the maximum solution of $S(a_i, b_i)$. We define $\overline{X} = \min_{i \in I}\{\hat{x}_i\}$.

**Definition 3.** Let $e: I \to J_i$ so that $e(i) = j \in J_i$, $\forall i \in I$, and let $E$ be the set of all vectors $e$. For the sake of convenience, we represent each $e \in E$ as an $m$-dimensional vector $e = [j_1, j_2, ..., j_m]$ in which $j_k = e(k)$. Also, we define $\underline{X}(e) = [\underline{X}(e)_1, \underline{X}(e)_2, ..., \underline{X}(e)_n] \in [0,1]^n$, where $\underline{X}(e)_j = \max_{i \in I}\{\breve{x}_i(e(i))_j\} = \max_{i \in I}\{\breve{x}_i(j_i)_j\}$, $\forall j \in J$.

Since Acze-Alsina t-norm is an operator with convex solution set (Lemma 1), we have the following theorem [5] that determines the feasible solutions set of problem (1).

**Theorem 1.** $S(A, b) = \bigcup_{e \in E} [\underline{X}(e), \overline{X}]$.

As a consequence, it turns out that $\overline{X}$ is the unique maximum solution and $\underline{X}(e)$'s ($e \in E$) are the minimal solutions of $S(A, b)$. Moreover, we have the following corollary that is directly resulted from Theorem 1.

**Corollary 1 (first necessary and sufficient condition).** $S(A, b) \neq \emptyset$ if and only if $\overline{X} \in S(A, b)$.

### 3. Linear programming subjected to Aczel-Alsina FRE

According to the well-known schemes used for optimization of linear problems such as (1) [5,10], problem (1) is converted to the following two sub-problems:

$$\begin{aligned} \min \ & Z_1 = \sum_{j=1}^n c_j^+ x_j & \min \ & Z_2 = \sum_{j=1}^n c_j^- x_j \\ & \varphi(a_i, x) = b_i \ , \ i \in I & & \varphi(a_i, x) = b_i \ , \ i \in I \\ & x \in [0,1]^n & & x \in [0,1]^n \end{aligned}$$

where $c_j^+ = \max\{c_j, 0\}$ and $c_j^- = \min\{c_j, 0\}$ for $j = 1, 2, ..., n$. It is easy to prove that $\overline{X}$ is the optimal solution of $Z_2$, and the optimal solution of $Z_1$ is $\underline{X}(e')$ for some $e' \in E$.

**Theorem 2.** Suppose that $S(A,b) \neq \emptyset$, and $\overline{X}$ and $\underline{X}(e^*)$ are the optimal solutions of sub-problems (3) and (2), respectively. Then $c^T x^*$ is the lower bound of the optimal objective function in (1), where $x^* = [x_1^*, x_2^*, ..., x_n^*]$ is defined as follows:

$$x_j^* = \begin{cases} \overline{X}_j & c_j < 0 \\ \underline{X}(e^*)_j & c_j \geq 0 \end{cases} \tag{2}$$

for $j = 1, 2, ..., n$.

**Proof.** Let $x \in S(A,b)$. Then, from Theorem 1 we have $x \in \bigcup_{e \in E}[\underline{X}(e), \overline{X}]$. Therefore, for each $j \in J$ such that $c_j \geq 0$, inequality $x_j^* \leq x_j$ implies $c_j^+ x_j^* \leq c_j^+ x_j$. In addition, for each $j \in J$ such that $c_j < 0$, inequality $x_j^* \geq x_j$ implies $c_j^- x_j^* \leq c_j^- x_j$. Hence, $\sum_{j=1}^{n} c_j x_j^* \leq \sum_{j=1}^{n} c_j x_j$. □

**Corollary 2.** Suppose that $S(A,b) \neq \emptyset$. Then, $x^* = [x_1^*, x_2^*, ..., x_n^*]$ as defined in (3), is the optimal solution of problem (1).

**Proof.** According to the definition of vector $x^*$, we have $\underline{X}(e^*)_j \leq x_j^* \leq \overline{X}_j$, $\forall j \in J$, which implies $x^* \in \bigcup_{e \in E}[\underline{X}(e), \overline{X}] = S(A,b)$. □

We now summarize the preceding discussion as an algorithm.

**Algorithm 1 (optimization of problem (1))**

Given problem (1):

**1.** Compute $J_i = \{j \in J : a_{ij} \geq b_i\}$ for each $i \in I$.

**2.** Compute $\overline{X} = \min_{i \in I}\{\hat{x}_i\}$ (Definition 2).

**3.** If $\overline{X} \notin S(A,b)$, then stop; $S(A,b)$ is empty (Corollary 1).

**4.** Find solutions $\underline{X}(e), \forall e \in E$ (Definition 3).

**5.** Find the optimal solution $\underline{X}(e^*)$ for $Z_1$ by pairwise comparison between the obtained solutions $\underline{X}(e)$.

**6.** Find the optimal solution $x^*$ for the problem (1) by (2) (Corollary 2).

## 4. Numerical examples

**Example 1.** Consider the following linear optimization problem (1):

$$\min \ Z = -7.6648 x_1 + 4.9208 x_2 + 6.1958 x_3 + 4.9047 x_4 - 3.2571 x_5$$
$$+ 1.6865 x_6 - 0.6209 x_7 - 8.2547 x_8$$

$$\begin{bmatrix} 0.1347 & 0.0955 & 0.0716 & 0 & 0.8463 & 0.0162 & 0.0115 & 0.1236 \\ 0.4505 & 0.1091 & 0.2857 & 0.4505 & 0.8606 & 0.4425 & 0.3448 & 0.6419 \\ 0 & 0.5548 & 0.0081 & 0.5723 & 0.9391 & 0.6595 & 0.6430 & 0.9200 \\ 0.7920 & 0.8793 & 0.7979 & 0 & 0.6802 & 0.2948 & 0.4479 & 0.3001 \\ 0.4197 & 0.2656 & 0.1975 & 0.4197 & 0.9174 & 0.9741 & 0.2847 & 0.0400 \\ 0.5325 & 0.8505 & 0.4725 & 0.5325 & 0.2567 & 0.9504 & 0.0982 & 0.7674 \end{bmatrix} \varphi x = \begin{bmatrix} 0.1347 \\ 0.4505 \\ 0.5723 \\ 0.792 \\ 0.4197 \\ 0.5325 \end{bmatrix}$$

$x \in [0,1]^8$

**Step 1:** In this example, we have $J_1 = \{1,5\}$, $J_2 = \{1,4,5,8\}$, $J_3 = \{4,5,6,7,8\}$, $J_4 = \{1,2,3\}$, $J_5 = \{1,4,5,6\}$ and $J_6 = \{1,2,4,6,8\}$.

**Step 2:** By Definition 1, we have $\hat{x}_1 = [1,1,1,1,0.1347,1,1,1]$, $\hat{x}_2 = [1,1,1,1,0.4513,1,1,0.4729]$, $\hat{x}_3 = [1,1,1,1,0.5724,0.6270,0.6413,0.5726]$, $\hat{x}_4 = [1,0.8030,0.8999,1,1,1,1,1]$, $\hat{x}_5 = [1,1,1,1,0.4198,0.4197,1,1]$ and $\hat{x}_6 = [1,0.5344,1,1,1,0.5325,1,0.5411]$. From Definition 2, $\overline{X} = [1, 0.5344, 0.8999, 1, 0.1347, 0.4197, 0.6413, 0.4729]$.

**Step 3:** Since $\overline{X} \in S(A,b)$, set $S(A,b)$ is feasible.

**Steps 4 and 5:** For this example $|E| = 2400$, that is, this problem has 2400 solutions $\underline{X}(e)$. By pairwise comparison, optimal solution $\underline{X}(e^*)$ for $Z_1$ is obtained as follows:

$$e^* = [1,1,7,1,1,1] \Rightarrow \underline{X}(e^*) = [1,0,0,0,0,0,0.6413,0]$$

**Step 6:** The optimal solution of Problem (1) is resulted as $x^* = [1,0,0,0,0.1347,0,0.6413,0.4729]$ with optimal objective value $Z^* = -12.4057$.

## Conclusion

Aczel-Alsina t-norm is an example of strict t-norms that are the most applied fuzzy operators in various fuzzy modelling problems. Particularly, Aczel – Alsina t-norm covers the most famous t-norms such as Minimum, Drastic and Product as $\lambda \to \infty$, $\lambda \to 0$ and $\lambda = 1$, respectively. In this paper, an algorithm was proposed to solve a linear optimization problem constrained with a special system of fuzzy relational equations (FRE) in which the feasible region was formed by the FRE defined by Aczel – Alsina t-norm. Based on the structural properties of the fuzzy equations, the feasible region of the problem is completely resolved and an algorithm was presented to find all the optimal solutions of the problem.